\documentclass[extra,referee]{gji}
\usepackage{timet,color}
\usepackage[urlcolor=blue,citecolor=black,linkcolor=black]{hyperref}
\usepackage{physics}
\usepackage{graphicx} 
\usepackage{algorithm}
\usepackage{algpseudocode}
\usepackage{amsmath} 
\usepackage{amssymb}
\usepackage{booktabs}
\usepackage{multirow}
\usepackage{microtype}

\newcommand{\mA}{\mathbf A}

\newcommand{\mJ}{\mathbf J}
\newcommand{\mK}{\mathbf K}
\newcommand{\mM}{\mathbf M}

\newcommand{\mQ}{\mathbf Q}

\newcommand{\bb}{\mathbf b}
\newcommand{\vf}{\mathbf f}
\newcommand{\vg}{\mathbf g}
\newcommand{\bu}{\mathbf u}
\newcommand{\ve}{\mathbf e}
\newcommand{\vm}{\mathbf m}
\newcommand{\vr}{\mathbf r}

\newcommand{\vv}{\mathbf v}

\title[Scalable Parallel 3-D TEM Inversion]{Scalable Parallel 3-D TEM Inversion via Rational Approximation of the Matrix Exponential}

\author[R.-U. Börner, S. Güttel, Th. Günther]
{Ralph-Uwe Börner$^1$, Stefan Güttel$^2$, Thomas Günther$^1$ \\
$^1$ Institute of Geophysics and Geoinformatics, TU Bergakademie Freiberg \\ Gustav-Zeuner-Straße 12, D-09599 Freiberg, Germany \\
$^2$ Department of Mathematics, The University of Manchester \\ Oxford Road, Manchester, M13 9PL, United Kingdom}

\date{ }

\usepackage{setspace}
\begin{document}

\maketitle

\begin{abstract}
We present a novel parallel implementation for large-scale three-dimensional electromagnetic inversion based on a Gauss–Newton framework combined with a rational near-best approximation of the matrix exponential for transient simulations. 
The method employs parallel direct solvers for the shifted linear systems arising from the partial fraction representation of the rational approximation and demonstrates efficient parallel execution on a shared-memory architecture using MPI. 
A key property of the approach is that the time dependence is entirely contained in the residuals of the employed rational functions, such that the computation of forward responses and sensitivities becomes effectively independent of the number of desired observation times. 
Model regularization is done with smoothness constraints, formulated with Raviart–Thomas $RT_0$ elements. 
The linearized inverse problems are solved using LSQR, using an implicit parallel Jacobian operator.
Numerical experiments demonstrate the successful recovery of a synthetic 3-D conductivity structure with approximately $700,\!000$ degrees of freedom. 
The study further discusses computational bottlenecks related to memory consumption and shared-memory scalability arising from the simultaneous storage of multiple sparse matrix factorizations. 
Possible improvements based on preconditioned iterative solvers and distributed high-performance computing architectures are outlined. 
The implementation in the Julia programming language is released as open-source software to support reproducible research and further development by the geophysical inversion community.
\end{abstract}


\section{Introduction}

Transient electromagnetics (TEM) is a well-established technique in applied geophysics, with applications in mineral exploration, hydrogeology and environmental studies.
In non-trivial environments with large-scale conductivity models, the inverse problem is inherently three-dimensional and therefore computationally demanding, as it requires repeated solutions of Maxwell’s equations in time domain.
Over the past two decades, considerable progress has been made in both forward modelling \citep[e.g., ][]{boerner2025rba,rochlitz2021evaluation} and inversion \citep[e.g., ][]{Haber2007,schwarzbach2013finite} methodologies, driven by advances in numerical discretization, solver technology and high-performance computing.

Current approaches to 3-D TEM inversion are mostly based on finite-element or finite-volume discretizations of the governing equations, combined with implicit time-integration schemes.

Within these frameworks, the transient response is computed by advancing the solution in time, typically using backward Euler type methods or related implicit schemes \cite[e.g.,][]{liu20163d}.
The inverse problem is commonly formulated as a regularized least-squares problem and solved using gradient-based (NLCG) or Hessian-based methods, such as Gauss--Newton or quasi-Newton algorithms, where sensitivities are obtained via adjoint-state techniques.
Finite-element formulations on unstructured meshes have proven particularly useful in this context, as they allow for the representation of complex geometries and support the construction of regularization operators consistent with the discretization \citep[e.g.,][]{Guenther2006} and furthermore allow for goal oriented refinement \citep{Spitzer2023}.

A number of practical implementations have demonstrated the feasibility of large-scale 3-D TEM inversion using these ingredients \citep[e.g.,][]{liu2019,liu2024three}.
\citet{liu20163d} proposed a parallel 3-D inversion scheme by formulating the forward problem in the frequency domain followed by a Fourier transform to obtain the EM response and the convolution with the transmitter waveform in the time domain.

Efficient use of sparse direct solvers and the reuse of matrix factorizations across time steps and sources have significantly reduced computational costs \citep{cai2022effective}.
Parallelization strategies typically exploit the independence of multiple sources \citep{liu2024three} and, to some extent, the reuse of system matrices within implicit time stepping.
By reconstructing the large sparse matrices arising in implicit time-stepping schemes into a series of low-order tridiagonal matrices, \citet{liu2025fast} accelerated the computation of both forward responses and the sensitivities.

Nevertheless, the computational workflow in these methods remains fundamentally organized around temporal marching, and the overall cost of both forward modelling and inversion is closely tied to the number of subsequently computed time steps required to represent the transient response.

The design of numerical methods for transient electromagnetic modelling is increasingly shaped by the limitations of Moore’s Law. 
As gains in single-core performance have plateaued, computational efficiency must now be achieved through parallelism rather than faster sequential execution. 
This development calls for a paradigm shift in algorithmic design: instead of adapting inherently sequential methods to parallel hardware, the mathematical formulation itself should expose concurrency.

As alternative viewpoint we consider the transient solution in terms of the action of a matrix exponential associated with the semi-discrete system. 
In this formulation, the time dependence of the electromagnetic fields is represented through an operator exponential. 
Recently, \citet{boerner2025rba} proposed a forward modelling approach based on uniform-in-time rational best approximations (RBA) of this matrix exponential. 
Their method constructs a family of rational functions sharing a common denominator, which, via partial fraction decomposition, gives rise to a set of shifted linear systems that can be solved independently.

This reformulation exemplifies the aforementioned paradigm shift. 
A key property is that the number of shifted systems is independent of the number of evaluation times, so that the computational cost is largely decoupled from the number of time channels. 
At the same time, the independent structure of these systems leads to an intrinsically parallel algorithm, in contrast to classical time-stepping or Krylov-based approaches, whose sequential nature limits scalability.

The inversion framework considered in this work is not new in itself, but rather serves as a template to demonstrate how such a reformulation of the forward problem propagates through the entire workflow. In particular, it illustrates how embedding parallelism at the level of the mathematical model enables scalable implementations of otherwise standard inversion techniques.

The proposed approach is not tied to a specific spatial discretization or regularization strategy.
Rather, it provides an alternative computational framework that can be combined with established finite-element formulations and standard regularization techniques.
In this sense, it complements existing inversion methodologies by introducing a different organization of the underlying computations, while retaining the flexibility required for realistic applications.
Our approach intends to enhance the performance of the codes available to the (time-domain) EM community.
As \citet{werthmueller2021towards} summarizes for CSEM, there is only a limited number of open-source 3-D modelling and inversion codes available like in the \textit{jInv} \citep{ruthotto2017jinv}, \textit{SimPEG} \citep{Cockett2015,Heagy2017}, \textit{custEM} \citep{rochlitz2019G,rochlitz2021evaluation}, or \textit{jMT3DAni} \citep{han2018julia} projects. 

Since the primary contribution of this work is methodological, we focus the numerical evalu\-ation on controlled synthetic experiments. 
Such experiments provide a fully specified setting in which the true model is known, enabling a direct and quantitative assessment of reconstruction accuracy, convergence behaviour, and computational performance. 
In addition, synthetic studies allow individual aspects of the inversion algorithm to be examined in isolation and under reproducible conditions, without the additional ambiguities introduced by field data, such as uncertain noise characteristics, imperfect system calibration, or survey geometry effects. 
As a result, they provide an effective framework for validating the underlying numerical method itself and for analysing its behaviour across a range of controlled scenarios. 
The results presented below demonstrate that the proposed approach reliably recovers the target structures while exhibiting the expected scaling and convergence properties.

An application to field data, while of practical interest, would primarily reflect additional modelling and acquisition uncertainties rather than provide further insight into the numerical properties of the proposed approach. Such investigations are therefore deferred to future work.

The remainder of this paper is organized as follows.
We first outline the mathematical formulation of the forward and inverse TEM problem and recall the rational approximation approach.
We then derive expressions for the sensitivities and describe the resulting inversion algorithm.
Numerical results for a synthetic 3-D model are presented to illustrate the performance of the method.
Finally, we summarize the main findings and discuss directions for further research.

\section{Theory}

\subsection{Governing equations}

We consider the quasi-static formulation of Maxwell’s equations for transient electromagnetic induction in a domain with a conductivity distribution. 
Eliminating the magnetic field leads to a second-order equation for the electric field $\ve = \ve(\vr, t)$ defined in a bounded domain $\Omega \subset \mathbb R^3$ such that
\begin{equation}
    \label{eq:pde}
    \nabla \times \mu^{-1} \nabla \times \ve + \sigma \partial_t \ve = -\partial_t \mathbf{j}^e, \quad\text{in } \Omega \times [0, \infty),
\end{equation}
where $\mu$ denotes magnetic permeability, $\sigma(\mathbf{r})$ the electrical conductivity, and $\mathbf{j}^e$ a prescribed source current density.
Appropriate boundary conditions are imposed on the outer boundary $\partial\Omega$ of the computational domain, typically enforcing a vanishing tangential electric field,
\begin{equation}
    \mathbf n \times \ve = \mathbf{0} \text{ along } \partial\Omega \times [0,\infty).
\end{equation}
As source term we consider a current density $\mathbf{j}^e$ resulting from a transmitter with a stationary current that is shut off at time $t = 0$,  resulting in 
\begin{equation}
    \mathbf{j}^e(\vr, t) = \mathbf q(\vr)H(-t),
\end{equation}
where $\mathbf q(\vr)$ denotes the spatial trace of the transmitter, and $H(\cdot)$ denotes the Heaviside step function.

In summary, the initial-value problem to be solved is  
\begin{subequations}
    \begin{align}
    \nabla \times \mu^{-1} \nabla \times \ve + \sigma \partial_t \ve  & = \mathbf{0} \quad \text{ on } \Omega \times (0, \infty) \label{eq:ivp} \\
    \mathbf{n} \times \ve & = \mathbf{0} \quad \text{ along } \partial\Omega \times (0, \infty) \label{eq:bc} \\
    \sigma \ve|_{t=0} & = \mathbf{q} \quad \text{ on } \Omega. \label{eq:source}
    \end{align}
\end{subequations}
This formulation forms the basis for both forward modelling and inversion.

\subsection{Spatial discretization}\label{sec:fe}

We discretize the operators in space using the finite-element method based on curl-conforming Nédélec elements on tetrahedral meshes \citep{schwarzbach2013finite}. 
The resulting semi-discrete system can be written in matrix form as
\begin{subequations}
\label{eq:ode}
\begin{align}
    \mK \bu(t) + \mM \partial_t \bu(t) & = \mathbf{0}, \quad t > 0 \\
    \mM \bu(0) & = \mathbf f.
\end{align}
\end{subequations}
Here, $\mathbf{u}(t) \in \mathbb R^N$ describes the $N$ degrees of freedom of the electric field, $\mathbf{K} \in \mathbb R^{N \times N}$ is the stiffness matrix associated with the curl–curl operator, and $\mathbf{M} \in \mathbb R^{N \times N}$ is the conductivity-weighted mass matrix.
Both matrices are sparse and inherit the properties of the underlying finite-element discretization.
Such formulations are standard in 3-D electromagnetic modelling and inversion; for details on their construction, properties and use in inverse problems we refer to the finite-element modelling frameworks on unstructured meshes \citep{schwarzbach2013finite,Heagy2017,rochlitz2021evaluation}.

The semi-discrete system can be rewritten as a linear system of ordinary differential equations
\begin{equation}
    \partial_t \bu (t) = -\mM^{-1} \mK \, \bu (t)
\end{equation}
whose formal solution is
\begin{equation}
    \bu(t) = \exp(-t \, \mM^{-1}\mK) (\mM^{-1}\vf) = \exp(-t \, \mM^{-1}\mK) \bb \quad\mbox{with}\quad\bb = \mM^{-1}\vf.
\label{eq:expm}
\end{equation}
This representation avoids explicit time discretization and expresses the transient response in terms of the action of a matrix exponential on a vector. 
In practice, the direct evaluation of this expression is not feasible for large-scale problems, and suitable approximations are required \citep{boerner2015arnoldi}.

Traditional approaches approximate the time evolution by implicit time stepping. In contrast, we adopt a rational approximation of the matrix exponential that is uniform in time.

\subsection{Rational approximation using RKFIT}

Following the approach of \cite{boerner2025rba}, the action of the matrix exponential is approximated by a family of rational functions constructed using the RKFIT (Rational Krylov Fitting) algorithm.

For a set of evaluation times $t_1,\dots,t_K$ we approximate $\exp(-t_j \mM^{-1} \mK) \bb$ by rational functions of the form
\begin{equation}
    \bu (t_j) \approx r_j(\mM^{-1} \mK) \bb = 2 \Re{\sum\limits_{i=1}^m \alpha_i^{[j]} (\mK - \xi_i \mM)^{-1} \vf }
\end{equation}
where $\{\xi_i\} \subset \mathbb C$ are the $m$ poles shared across all time points and $\{ \alpha_i^{[j]} \} \subset \mathbb C$ are  coefficients that are precomputed for the chosen times $t_j$ and poles $\xi_i$, independent on any spatial discretization.
As crucial feature of this formulation, the shifted system matrices $(\mathbf{K} - \xi_i \mathbf{M})$ are independent of time.
As a consequence, the factorization of these matrices can be reused for all evaluation times.

\begin{figure}
    \centering\includegraphics[width=\linewidth]{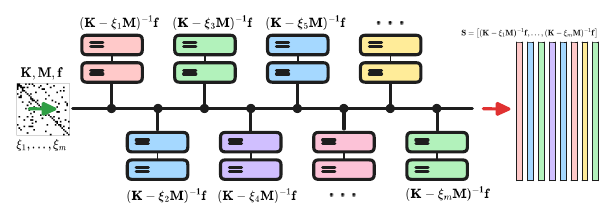}
    \caption{Schematic illustration of the parallel distribution of shifted linear systems arising from the rational approximation. The system matrices $\mathbf{K}$, $\mathbf{M}$, the poles ${\xi_i}$, and the right-hand side vector $\mathbf{f}$ are distributed across computational resources. Each worker is assigned one or more poles $\xi_i$, forming, factorizing, and solving the corresponding linear systems $(\mathbf{K} - \xi_i \mathbf{M})\mathbf{g}_i = \mathbf{f}$ independently. 
    The resulting solutions $\mathbf g_i$ are subsequently gathered on a master process, where the final field response is obtained by forming the linear combination with the corresponding rational approximation residuals.}\label{fig:distribution}
\end{figure}

This leads to two important properties: 
First, the number of required linear system to be solved depends only on the number of poles $m$, and not on the number of evaluation times $K$.
Second, the shifted systems are independent and can be solved in parallel, providing a naturally scalable implementation.
An overview of the parallel solution strategy is provided in Algorithm~\ref{alg:parallel-mumps}, and is illustrated in Fig.~\ref{fig:distribution}.

In the following, we exploit this representation not only for forward modelling but also for the computation of sensitivities and the solution of the inverse problem.
By expressing both the forward response and its derivatives in terms of the same set of shifted systems, the computational structure of the forward solver is retained within the inversion.

\begin{algorithm}
\caption{Parallel solution using cached MUMPS factorizations}
\label{alg:parallel-mumps}
\begin{algorithmic}[1]
\Require Matrices $\mK$, $\mM$, right-hand side $\vf$, poles $\{\xi_i\}_{i=1}^{m}$
\Ensure Local solutions $\vg_i = (\mK-\xi_i \mM)^{-1}\vf$ 

\State Determine MPI rank $p$ and local pole set $\mathcal I_p \subset \{1,\ldots,m\}$

\ForAll{$i \in \mathcal I_p$}
    \State Form the shifted matrix
    \[
        \mA_i \gets \mK - \xi_i \mM
    \]

    \If{factorization of $\mA_i$ is not cached}
        \State Factorize $\mA_i$ with MUMPS on the local rank
        \State Store the factorization in the local cache
    \EndIf

    \State Compute
    \[
        \mA_i^{-1} \vf =  \vg_i
    \]
    using the cached MUMPS factorization
    \State Store $\vg_i$ as the local solution associated with pole $\xi_i$
\EndFor

\State Gather or expose the distributed collection $\{\vg_i\}_{i=1}^{m}$
\State Compute time-dependent responses by linear combination:
\[
    u(t_j) = 2\operatorname{Re}
    \left(
        \sum_{i=1}^{m} \alpha_i^{[j]} \vg_i
    \right),
    \qquad j=1,\ldots,K .
\]

\end{algorithmic}
\end{algorithm}

\section{The Inverse Problem}

Let $\vm \in \mathbb R^P$ denote the model parameter vector, typically representing the natural logarithm of the electrical conductivity, $\mathbf{m} = \log \boldsymbol{\sigma}$, defined piecewise over the computational mesh.
For a given model $\mathbf{m}$, the semi-discrete system introduced in Section 2 yields a time-dependent solution 
$\mathbf{u}(t,\mathbf{m})$.

Measurements are related to the state variables through an observation operator $\mathbf{Q} \in \mathbb R^{M \times N}$, such that the predicted data at time $t_j$ are given by
\begin{equation}
\mathbf{d}_j(\mathbf{m}) = \mathbf{Q} \, \mathbf{u}(t_j,\mathbf{m}), \quad j = 1, \dots, K.
\end{equation}

Collecting the time channels, we write
\begin{equation}
\mathbf d (\vm) =
\begin{bmatrix}
     \mathbf d_1 (\vm) \\ \vdots \\ \mathbf d_K (\vm)
\end{bmatrix} \in \mathbb R^{M K}.
\end{equation}

Given observed data $\mathbf{d}^{\mathrm{obs}}$, the inverse problem is formulated as the minimization of a regularized least-squares functional
\begin{equation}
\phi(\mathbf{m}) = \frac{1}{2} \| \mathbf{W}_d \bigl(\mathbf{d}(\mathbf{m}) - \mathbf{d}^{\mathrm{obs}} \bigr) \|_2^2 + \lambda \, R(\mathbf{m}),
\end{equation}
where $\mathbf{W}_d$ is a data-weighting matrix, $R(\mathbf{m})$ is a regularization functional, and $\lambda > 0$ is a regularization parameter.
The choice of $R(\mathbf{m})$ is not specific to the present formulation and may follow standard approaches, such as smoothness or total variation regularization.

The minimization of $\phi(\mathbf{m})$ is typically carried out using a gradient-based method, for example a {Gauss-Newton} or quasi-Newton scheme, which requires the evaluation of the Jacobian of the data with respect to the model parameters.

\subsection{Jacobian computation}

To compute the Jacobian of the predicted data with respect to the model parameters, we differentiate the forward response with respect to $\mathbf{m}$. 
For each pole $\xi_i$, we define
\begin{equation}
\mathbf{g}_i = \mathbf{A}_i(\mathbf{m})^{-1} \mathbf{f} := [\mK - \xi_i \mM(\vm)]^{-1} \vf.
\end{equation}

Using the identity for the derivative of the inverse matrix,
\begin{equation}
\partial_{\mathbf{m}} \mathbf{A}_i^{-1} = - \mathbf{A}_i^{-1} \, (\partial_{\mathbf{m}} \mathbf{A}_i) \, \mathbf{A}_i^{-1},
\end{equation}
and noting that
\begin{equation}
\partial_{\mathbf{m}} \mathbf{A}_i(\mathbf{m}) = -\xi_i \, \partial_{\mathbf{m}} \mathbf{M}(\mathbf{m}),
\end{equation}
we obtain the sensitivity of the predicted data at time $t_j$ as
\begin{equation}
\mathbf{J}_j(\mathbf{m}) = \partial_{\mathbf{m}} \mathbf{d}_j(\mathbf{m})
= 2 \Re{\sum_{i=1}^{{m}} \alpha_i^{[j]} \, \xi_i \, \mathbf{Q} \, \mathbf{A}_i^{-1}
\bigl[ \partial_{\mathbf{m}} \mathbf{M}(\mathbf{m}) \, \mathbf{g}_i \bigr]}.
\end{equation}

The full Jacobian is obtained by stacking the contributions from all time channels.
We obtain
\begin{equation}
    \mathbf J := 
    \begin{pmatrix}
    \mathbf J_1 \\ \mathbf J_1 \\ \vdots \\ \mathbf J_K    
    \end{pmatrix}
    = 2 \Re \left(
    \sum_{i=1}^m 
    \begin{pmatrix}
        \alpha_i^{[1]} \\ \alpha_i^{[2]} \\ \vdots \\ \alpha_i^{[K]}
    \end{pmatrix}
    \otimes \xi_i \mathbf Q \mA_i(\vm)^{-1} [\partial_m \mM (\vm) \mathbf g_i]
    \right),
\end{equation}
where $\otimes$ denotes the Kronecker tensor product.

In practice, the Jacobian does not need to be formed explicitly, since both actions $\mathbf{J}\mathbf{v}$ and $\mathbf{J}^\top \mathbf{w}$ can be evaluated through the same set of shifted linear systems using the stored factorizations.
The action of $\mathbf J$ on the vector $\vv$ can be written as
\begin{equation}
    \mathbf J \vv := 
    \begin{pmatrix}
    \mathbf J_1 \vv \\ \mathbf J_1 \vv \\ \vdots \\ \mathbf J_K \vv    
    \end{pmatrix}
    = 2 \Re \left(
    \sum_{i=1}^m 
    \begin{pmatrix}
        \alpha_i^{[1]} \\ \alpha_i^{[2]} \\ \vdots \\ \alpha_i^{[K]}
    \end{pmatrix}
    \otimes \xi_i \mathbf Q \mA_i(\vm)^{-1} [\partial_m \mM (\vm) \mathbf g_i \vv]
    \right).
\end{equation}

The derivative $\partial_m \mM(\vm)$ is a very sparse third-order tensor of dimension $N \times N \times P$, where each slice (with respect to a model parameter) corresponds to the contribution of a single element mass matrix assembled with unit conductivity that is stored at the associated global indexing scheme of the degrees of freedom.

Contracting this tensor with a vector $\mathbf g_i$ amounts to a mode-2 tensor–vector product 
$(\partial_m \mM)\, \mathbf g_i$ of dimension $N \times P$, which combines the element-wise mass matrix contributions with the vector~$\mathbf g_i$. 
The result is a sparse rectangular matrix: each column corresponds to a model parameter (element), and each row reflects the action of the corresponding local mass matrix on $\mathbf g_i$. 
Its column dimension is therefore equal to the number of model parameters and matches the dimension required for subsequent multiplication with the vector $\mathbf v$ in the Jacobian expression.
The resulting term $\partial_m \mM (\vm) \mathbf g_i \vv$ is a complex vector of dimension ${N}$.

The adjoint Jacobian $\mathbf J_j^\top$ for the time channel $j$ is given by
\begin{equation}
    \mathbf J_j^\top := 2 \Re \left(
    \sum_{i=1}^m (\alpha_i^{[1]} \ \alpha_i^{[2]} \dots \ \alpha_i^{[K])}) \otimes 
    \xi_i [\partial_m \mM (\vm) \mathbf g_i ]^\top \mA_i^{-1} \mathbf Q^\top
    \right),
\end{equation}
which is of dimension $P \times M$.
A suitable arrangement of the data vector $\mathbf w$ is a partition of the form
\begin{equation}
    \mathbf w = \begin{pmatrix}
        \overline{\mathbf w}_1 \\
        \overline{\mathbf w}_2 \\
        \vdots \\
        \overline{\mathbf w}_K
    \end{pmatrix} \in \mathbb R^{M K},
\end{equation}
such that
\begin{equation}
    \mathbf J^\top \mathbf w = \sum_{j=1}^K \mathbf J_j^\top \overline{\mathbf w}_j \in \mathbb R^P.
\end{equation}

The correctness of these operators can be verified by Taylor remainder tests and by the standard adjoint test, as discussed below.
Let $\delta \mathbf{m}$ denote an arbitrary perturbation of the model parameters and let $h > 0$ be a scalar step length. 
If the Jacobian $\mathbf{J}(\mathbf{m})$ is implemented correctly, then the perturbed predicted data satisfy
\begin{equation}
\mathbf{d}(\mathbf{m} + h \, \delta \mathbf{m})
=
\mathbf{d}(\mathbf{m})
+
h \, \mathbf{J}(\mathbf{m}) \, \delta \mathbf{m}
+
\mathcal{O}(h^2)
\quad \text{as } h \to 0 .
\end{equation}

Consequently, the residual
\begin{equation}\label{eq:e1}
\mathbf{r}_1(h) = \mathbf{d}(\mathbf{m} + h \, \delta \mathbf{m})
- \mathbf{d}(\mathbf{m})
- h \, \mathbf{J}(\mathbf{m}) \, \delta \mathbf{m}
\end{equation}
should satisfy
\begin{equation}
e_1(h) := \| \mathbf{r}_1(h) \|_2 = \mathcal{O}(h^2).
\end{equation}

In practice, this behaviour is verified by evaluating the residual for a sequence of decreasing step lengths $h$ and confirming second-order decay over the range in which truncation and round-off errors remain negligible.

\begin{figure}
    \centering\includegraphics[width=0.7\linewidth]{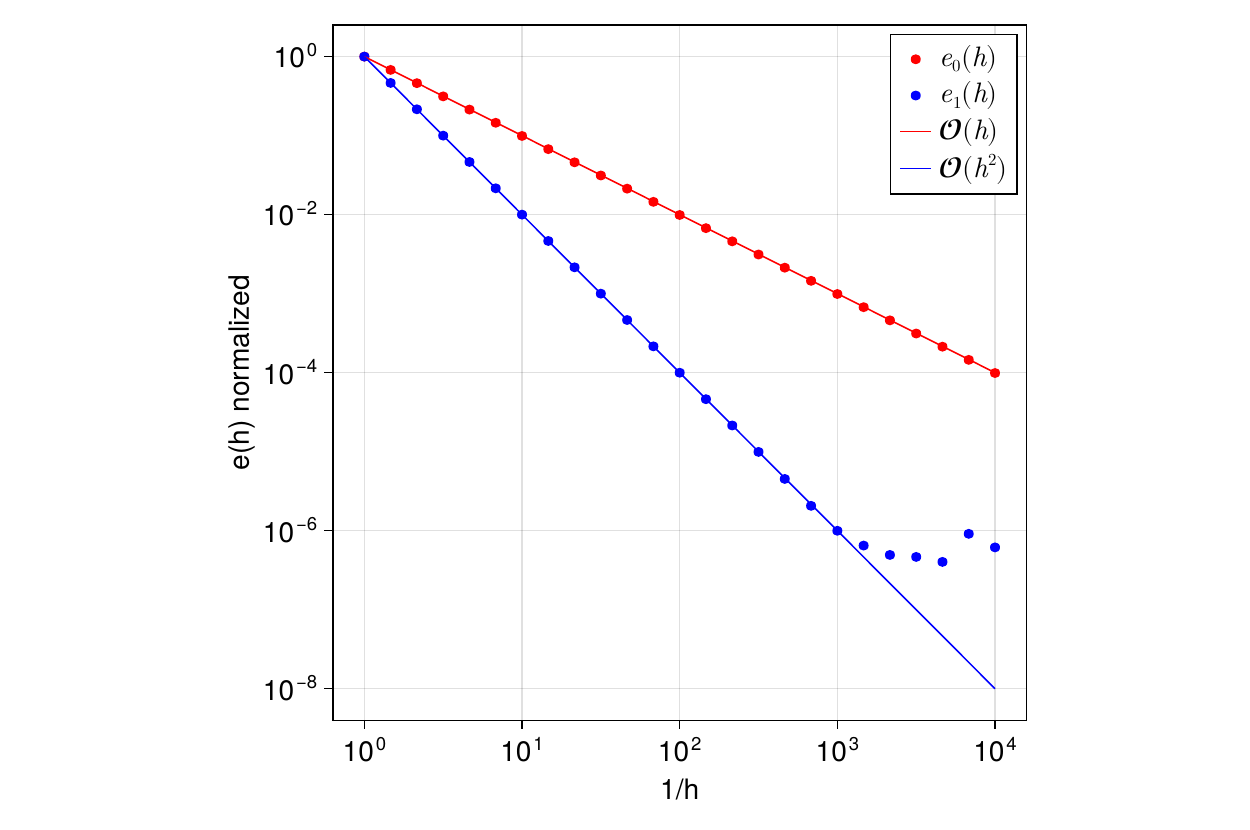}
    \caption{Taylor remainder test for the forward operator and its linearization. Red markers show the norm of the first-order difference $e_0(h)$ (\ref{eq:e0}), while the blue markers show the norm of the second-order remainder $e_1(h)$ (\ref{eq:e1}). Solid lines indicate the expected asymptotic behaviour of order $\mathcal{O}(h)$ and $\mathcal{O}(h^2)$, respectively.}
    \label{fig:taylortest}
\end{figure}

A complementary check is obtained from the first-order remainder
\begin{equation}\label{eq:e0}
\mathbf{r}_0(h)
=
\mathbf{d}(\mathbf{m} + h \, \delta \mathbf{m})
-
\mathbf{d}(\mathbf{m}),
\end{equation}
for which one expects
\begin{equation}\label{eq:e0}
e_0(h) := \| \mathbf{r}_0(h) \|_2 = \mathcal{O}(h).
\end{equation}
Thus, in a log-log representation, the uncorrected difference should exhibit a slope of $-1$, whereas the Jacobian-corrected remainder should exhibit a slope of $-2$ with respect to $h^{-1}$. 
This provides a straightforward numerical verification of the linearization as shown in Fig. \ref{fig:taylortest}.

In addition, the consistency of the adjoint implementation can be checked by the standard adjoint test \citep[p. 88]{haber2014book}.

Let $\delta \mathbf{m}$ be an arbitrary model perturbation and let $\delta \mathbf{d}$ be an arbitrary data-space vector of compatible dimension. 
Then the discrete Jacobian and its adjoint must satisfy
\begin{equation}
\langle \mathbf{J} \, \delta \mathbf{m}, \delta \mathbf{d} \rangle
=
\langle \delta \mathbf{m}, \mathbf{J}^\top \delta \mathbf{d} \rangle
\end{equation}
up to machine precision, where $\langle \cdot,\cdot \rangle$ denotes the Euclidean inner product. 
In numerical experiments, the relative mismatch
\begin{equation}
\frac{
\left|
\langle \mathbf{J} \, \delta \mathbf{m}, \delta \mathbf{d} \rangle
-
\langle \delta \mathbf{m}, \mathbf{J}^\top \delta \mathbf{d} \rangle
\right|
}{
\max\!\left(
|\langle \mathbf{J} \, \delta \mathbf{m}, \delta \mathbf{d} \rangle|,
|\langle \delta \mathbf{m}, \mathbf{J}^\top \delta \mathbf{d} \rangle|
\right)
}
\end{equation}
should be close to machine accuracy.

These two tests provide complementary validation of the sensitivity implementation.
It should be noted that these tests establish the algebraic consistency of the linearization rather than the physical correctness of the model. Any systematic errors in the forward operator are inherited by the Jacobian, such that the results confirm the internal consistency of the discrete formulation.
For example, insufficient spatial resolution in the vicinity of receiver locations may lead to inaccuracies in the evaluation of the observation operator, resulting in model responses that are affected by discretization artefacts. 
Such effects are not detected by the Taylor or adjoint tests, as they are consistently propagated through both the forward operator and its Jacobian.

A key feature of the Jacobian calculation is that the vectors $\mathbf{g}_i$ depend only on the poles $\xi_i$ and not on the time index $j$. 
Consequently, the computation of sensitivities requires:
\begin{enumerate}
    \item the parallel solution of $m$ forward problems to obtain $\mathbf{g}_i$,
	\item the application of the derivative of the mass matrix to these vectors, and
	\item additional parallel solves with the same matrices $\mathbf{A}_i(\mathbf{m})$ based on their respective factorizations.
\end{enumerate}

All linear systems involved correspond to the same set of shifted matrices $\mathbf{A}_i(\mathbf{m})$ as in the forward problem.
The formulation above reveals that both the forward response and its sensitivities can be expressed in terms of a fixed set of shifted linear systems associated with the poles $\{\xi_i\}$. This has two important consequences:
\begin{enumerate}
    \item the number of required system solutions is independent of the number of time channels $K$. 
Once the solutions $\mathbf{g}_i$ have been computed, the response at all times is obtained through linear combinations with coefficients $\alpha_i^{[j]}$;
    \item the shifted systems are independent and can be solved in parallel. 
The same factorizations can be reused for both forward modelling and sensitivity calculations, which leads to a consistent computational framework for inversion.
\end{enumerate}

In contrast to time-stepping-based approaches, where sensitivities are typically computed by forward and adjoint solutions at each time step, the present formulation avoids temporal recursion. 
Instead, the inversion is organized around the evaluation of a small number of model-dependent linear systems, which constitutes the dominant computational cost.

\subsection{Regularization}

The inverse problem is ill-posed and requires regularization to obtain stable and physically meaningful solutions. 
In this work, we employ a smoothness regularization following the finite-element formulation introduced by \cite{schwarzbach2013finite}.

We consider the model parameter $\mathbf{m} = \log \boldsymbol{\sigma}$ represented as a piecewise constant function on the elements of the mesh. 
In this setting, the classical gradient $\nabla \mathbf{m}$ is not defined in the standard sense, as discontinuities occur across element interfaces. 
To address this, a mixed finite-element formulation is employed in which the gradient is represented in a weak sense using a dual variable.

Specifically, the gradient of the model is approximated in the space $H(\mathrm{div})$ using lowest-order Raviart--Thomas ($RT_0$) elements. 

The smoothness functional can then be written in discrete form as
\begin{equation}
R(\mathbf{m}) = \frac{1}{2} (\mathbf{m} - \mathbf{m}^{\mathrm{ref}})^\top \mathbf{L} (\mathbf{m} - \mathbf{m}^{\mathrm{ref}}),
\end{equation}
where the regularization matrix $\mathbf{L}$ is given by
\begin{equation}
\mathbf{L} = \mathbf{D} \, \mathbf{M}_{\mathrm{div}}^{-1} \, \mathbf{D}^\top.
\end{equation}
Here, $\mathbf{D}$ is a discrete divergence operator mapping from face-based $RT_0$ degrees of freedom to cell-based quantities, and $\mathbf{M}_{\mathrm{div}}$ is the corresponding mass matrix in the $H(\mathrm{div})$ space.

This construction provides a consistent discretization of the model gradient for piecewise constant parameters on unstructured meshes and avoids the ambiguities associated with defining finite differences across irregular element interfaces. 
In particular, the resulting regularization operator does not require heuristic scaling with element size or face area.

The gradient and Hessian of the regularization functional are given by
\begin{equation}
\nabla_m R(\mathbf{m}) = \mathbf{L} (\mathbf{m} - \mathbf{m}^{\mathrm{ref}}), \quad
\nabla^2_m R(\mathbf{m}) = \mathbf{L}.
\end{equation}

In practice, the inversion of the mass matrix $\mathbf{M}_{\mathrm{div}}$ is approximated by a diagonal or lumped representation to reduce computational cost, without significantly affecting the qualitative behaviour of the regularization.

\subsection{Optimization}

The minimization of the objective functional is carried out within a Gauss–Newton framework. At iteration $\nu$, the model update $\delta \mathbf{m}$ is obtained from the linearized least-squares problem
\begin{equation}\label{eq:lsqr}
\min_{\delta \mathbf{m}}
\;
\left\|
\begin{bmatrix}
\mathbf{W}_d\mathbf{J} \\
\sqrt{\lambda}\,\mathbf{L}^{1/2}
\end{bmatrix}
\delta \mathbf{m}
+
\begin{bmatrix}
\mathbf{W}_d\bigl(\mathbf{d}(\mathbf{m}^{(\nu)})-\mathbf{d}^{\mathrm{obs}}\bigr) \\
\sqrt{\lambda}\,\mathbf{L}^{1/2}\bigl(\mathbf{m}^{(\nu)}-\mathbf{m}^{\mathrm{ref}}\bigr)
\end{bmatrix}
\right\|_2^2 ,
\end{equation}
where $\mathbf{J}$ denotes the Jacobian evaluated at the current iterate $\mathbf{m}^{(\nu)}$, $\mathbf{W}_d$ is the data-weighting matrix, and $\mathbf{L}$ is the regularization operator introduced above.
Rather than forming the normal equations 
\begin{equation}\label{eq:cgls}
(\mathbf{J}^\top \mathbf{W}_d^\top \mathbf{W}_d \mathbf{J} + \lambda \mathbf L) \, \delta \mathbf{m} = 
\mathbf{J}^\top \mathbf{W}_d^\top \mathbf{W}_d (\mathbf{d}^{\mathrm{obs}}-\mathbf{d}(\mathbf{m}^{(\nu)}))
+ \lambda \mathbf{L} (\mathbf{m}^{\mathrm{ref}} - \mathbf{m}^{(\nu)}),
\end{equation}
we solve the least-squares problem (\ref{eq:lsqr}) using the LSQR algorithm \citep{paige1982lsqr}.
This is advantageous, particularly for parallel computation \citep{Lee2013}, compared to solving (\ref{eq:cgls}) with conjugate gradients (CGLS) as done by \citet{Guenther2006} as the condition number is effectively squared relative to that of the original operator $\mathbf{W}_d\mathbf{J}$.
This amplification of ill-conditioning can lead to slower convergence and increased sensitivity to numerical rounding errors which is avoided by LSQR.

In our implementation, the Jacobian is represented in matrix-free form through a linear operator. More precisely, for a given vector $\mathbf{v}$ in model space, the product
$\mathbf{J}\mathbf{v}$ is evaluated by a function that applies the action of the Jacobian on a vector using the rational approximation framework, which amounts to solving the corresponding shifted linear systems. 
Similarly, for a given vector $\mathbf{w}$ in data space, the adjoint action
$\mathbf{J}^\top \mathbf{w}$
is computed through the corresponding adjoint solves. 
These two operator evaluations are provided to LSQR through a linear map abstraction, such that the iterative solver treats the Jacobian as an implicit operator rather than an assembled matrix.

The application of the regularization operator is simplified by the use of a lumped approximation of the mass matrix $\mathbf{M}_{\mathrm{div}}$ in the $RT_0$ formulation. 
This results in a symmetric positive definite regularization matrix $\mathbf{L}$, for which a Cholesky factorization
$\mathbf{L} = \mathbf{R}^\top \mathbf{R}$
can be computed. 
The action of $\mathbf{L}^{1/2}$ required in the augmented system is then realized through multiplication with the factor~$\mathbf{R}$.
This avoids the explicit and expensive construction of matrix square roots and allows for an efficient evaluation of the regularization term within the LSQR iterations.

The Gauss--Newton method provides a descent direction $\delta \mathbf{m}$, but a full update step may not always lead to a decrease of the objective functional. 
To ensure robust convergence, we employ a backtracking line search. Starting from an initial step length $\eta = 1$, the trial update
\begin{equation}
\mathbf{m}^{(\nu+1)} = \mathbf{m}^{(\nu)} + \eta \, \delta \mathbf{m}
\end{equation}
is successively halved until the \textit{Armijo} condition
\begin{equation}
\phi(\mathbf{m}^{(\nu)} + \eta \, \delta \mathbf{m})
\le
\phi(\mathbf{m}^{(\nu)}) + c_1 \eta \, \nabla \phi(\mathbf{m}^{(\nu)})^\top \delta \mathbf{m}
\end{equation}
is met, using a constant $c_1$ of $10^{-4}$ \citep{nocedal2006numerical}.

In practice, the step length $\eta$ in the line search must be chosen with care.
Although the exact Gauss–Newton direction is a descent direction, this guarantee does not extend to the inexact setting considered here, where the linearized system is solved approximately using LSQR and includes the regularization term, weighted by $\lambda$. 
The value of $\lambda$ directly affects the conditioning of the system and the relative scaling of the update components.

We observe that for excessively small $\eta$, the interplay between Krylov solver inaccuracies and the regularization scaling may lead to a loss of descent, with the update failing to reduce the objective function.
Accordingly, the line search is safeguarded by enforcing a lower bound on $\eta$.

Since the gradient $\nabla_m \phi(\mathbf{m}^{(\nu)})$ 
at $\eta=0$ is available at no additional cost, the evaluation of $\phi$ at a trial step $0 < \eta < 1$ also allows for a local quadratic model of the objective function along the search direction. 
The minimizer of this quadratic can be used as an estimate for an improved step length. 
In the present implementation, this possibility is used only as a heuristic, while the final step length is determined by the \textit{Armijo} criterion.

The regularization strength $\lambda$ controls the smoothness of the model and needs to be chosen such that the data are fitted within error bounds, but also affect the convergence of the inverse problem.
We follow an adaptive cooling strategy with piece-wise constant values:
Starting with an initial value, the objective function is minimized until convergence is observed.
Unless the data fit approaches $\chi^2\approx 1$, $\lambda$ is halved and the inversion is continued.

\section{Synthetic examples}

\subsection{Computational setup}

All numerical experiments were carried out using a finite-element implementation based on the Julia programming language \citep{bezanson2017julia} . 
The spatial discretization employs curl-conforming Nédélec elements as described in Section \ref{sec:fe}, implemented using the \textit{Gridap.jl} framework \citep{badia2020gridap}, which provides a flexible environment for the assembly of the discrete operators and the formulation of the forward problem.
The computational meshes were generated with the mesh generator \textit{Gmsh} \citep{geuzaine2009gmsh}, allowing for flexible irregular tetrahedral discretizations adapted to the geometry of any topography, anomalies and the transmitter–receiver configuration.
Most of the figures were prepared using \textit{Makie.jl}, a flexible high-performance data visualization package for Julia \citep{danisch2021makie}. 

The solution of the shifted linear systems arising from the rational approximation is performed using sparse direct solvers such as  the widespread \citep[cf. overview in][]{weiss2025evaluation} sparse multi-frontal solver \textit{MUMPS} \citep{MUMPS:1}.
To exploit the inherent parallelism of the formulation, we have developed a custom Julia interface to \textit{MUMPS}, implemented using the Message Passing Interface \textit{MPI}. 
This interface enables the concurrent factorization and solution of multiple shifted systems corresponding to the poles of the rational approximation.
All computations were performed on a shared-memory system equipped with four Intel Xeon processors, providing 48 physical cores and 48 virtual cores, and a total of 2 TB of RAM. 

The independent shifted systems were distributed across processes, which were pinned to individual CPU sockets to avoid contention and to ensure reproducible performance. 
This strategy allows the parallel structure of the algorithm to be exploited without relying on fine-grained threading within individual solver instances.

A factorization of one shifted system typically exceeds several GB of memory. As a consequence, data access is often limited by memory bandwidth rather than compute performance, and cache locality is reduced. This may lead to a decrease in parallel efficiency, particularly when multiple factorizations are performed concurrently.

In addition, the shared-memory architecture used in our experiments exhibits non-uniform memory access (NUMA) characteristics. 
Memory attached to a given CPU socket can be accessed with lower latency and higher bandwidth by processes running on the same socket than by processes on other sockets. 
Consequently, memory access across sockets (e.g., when a process on one CPU accesses data residing in memory attached to a CPU on another socket) incurs additional latency and reduced bandwidth. 
For the large factorizations considered here, which exceed typical cache sizes by several orders of magnitude, such NUMA effects can lead to further performance degradation if memory locality is not preserved. 
The explicit pinning of processes to sockets partially mitigates this effect, but does not eliminate it entirely when data are distributed across the global memory space.

Given the independence of the shifted systems, further performance gains can be expected when distributing these computations across nodes on a high-performance cluster.
In such settings, the cache over-subscription per node is reduced and communication between tasks is minimal, suggesting favourable scalability for large-scale problems.
Our implementation reflects the computational structure of the proposed method, in which the dominant cost arises from a fixed number of independent linear system solves that can be executed concurrently.

\subsection{Synthetic model description}

\begin{figure}
    \centering
    \includegraphics[width=1.0\linewidth]{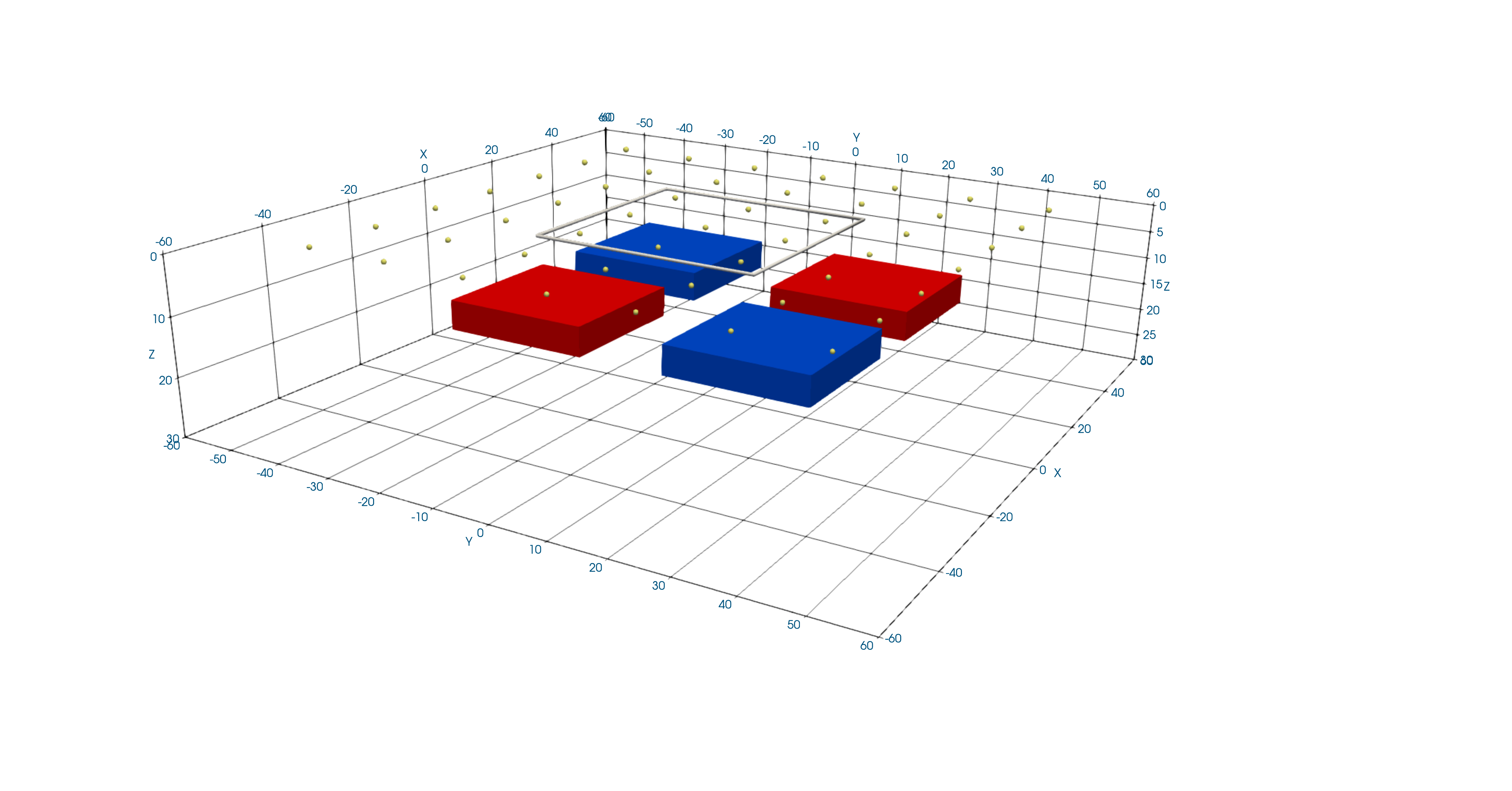}
    \vspace*{-12mm}
    \caption{Three-dimensional view of the synthetic model and survey configuration. The square transmitter loop is located at the surface, together with the receiver grid, while two conductive (red) and two resistive (blue) blocks are embedded in a homogeneous half-space at shallow depth. The uniform background is omitted for clarity.}
    \label{fig:model-view}
\end{figure}

To assess the performance of the proposed inversion framework, we consider a synthetic model comprising a homogeneous half-space with an electrical conductivity of 0.1~$\mathrm{S/m}$, into which two conductive and two resistive anomalies are embedded.
Each anomaly is represented by a rectangular block of dimensions 25 $\times$ 25 $\times$ 5~$\mathrm{m}^3$, with conductivities of 1~$\mathrm{S/m}$ and 0.01~$\mathrm{S/m}$, respectively.
The blocks are placed at a depth of 10 m below the surface (top of the bodies) and are laterally distributed within the computational domain to form a non-trivial three-dimensional target structure (Fig.~\ref{fig:model-view}).

The transmitter is a square loop of side length 40 m, located at the surface ($z=0$). 
The current is assumed to be switched off at $t=0$, and the transient response is computed for $t>0$. 
The receiver array consists of 49 surface stations arranged in a regular $7 \times 7$ grid with a spacing of 15\,m, covering an area from $-$60\,m to 60\,m in both horizontal directions (Fig.~\ref{fig:model-view}). 
In our numerical experiment, the number of degrees of freedom amounts to $N = 685{,}350$.

The measured quantity is the time derivative of the vertical magnetic flux density, $\partial_t B_z$, evaluated at the receiver locations.
For each receiver, the containing finite element is identified, and the observation operator $\mQ$ is assembled by evaluating the curl of the corresponding local basis functions at the receiver position and extracting the vertical component.
Applied to the discrete solution, this provides an approximation of $\partial_t B_z$.

At each receiver location, the transient response is obtained at 31 logarithmically spaced time channels in the interval $t \in [10^{-6}, 10^{-3}]$ s. 
This time range is representative for near-surface TEM surveys and is consistent with the interval used in the construction of the rational approximation.

Synthetic data are generated by evaluating the forward model for the true conductivity distribution. 
Gaussian random noise is added to the synthetic observations, using a standard deviation of an error model $\sigma(d) = |d|\,\epsilon_r + \epsilon_a$, with relative and absolute error levels $\epsilon_r$ and $\epsilon_a$, respectively.
The data are weighted accordingly by scaling with the inverse of $\sigma(d)$.

\subsection{Reference inversion example}

\begin{figure}
    \centering
    \includegraphics[width=\linewidth]{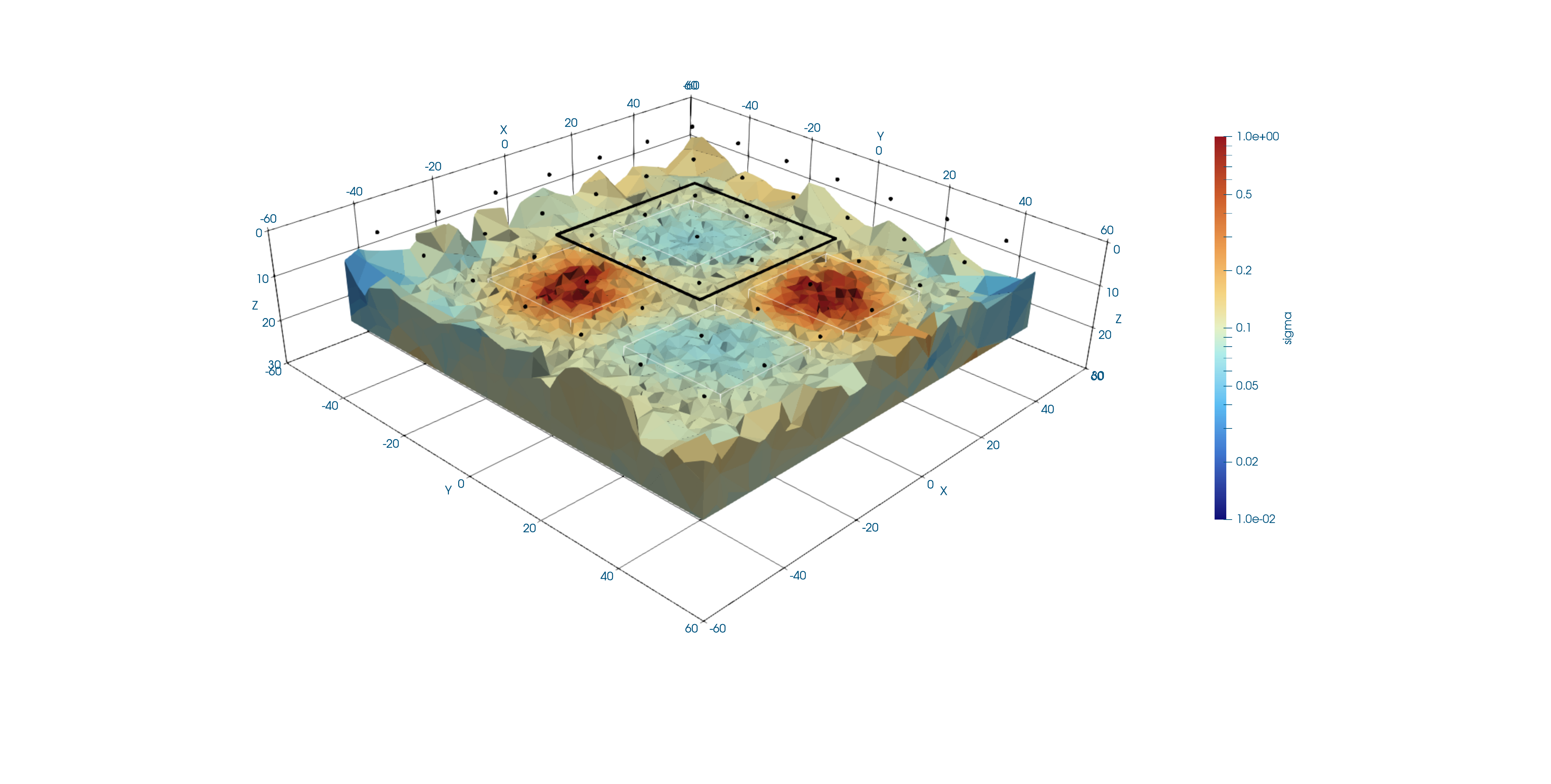}
    \vspace*{-12mm}
    \caption{Three-dimensional view of the reconstructed conductivity model obtained from the synthetic data after 25 iterations. The inversion recovers the locations and general shapes of the conductive blocks embedded in the resistive half-space. The anomalies appear slightly smoothed and reduced in amplitude due to regularization effects.}
    \label{fig:model-result}
\end{figure}

\begin{figure}
    \centering
    \includegraphics[width=\linewidth]{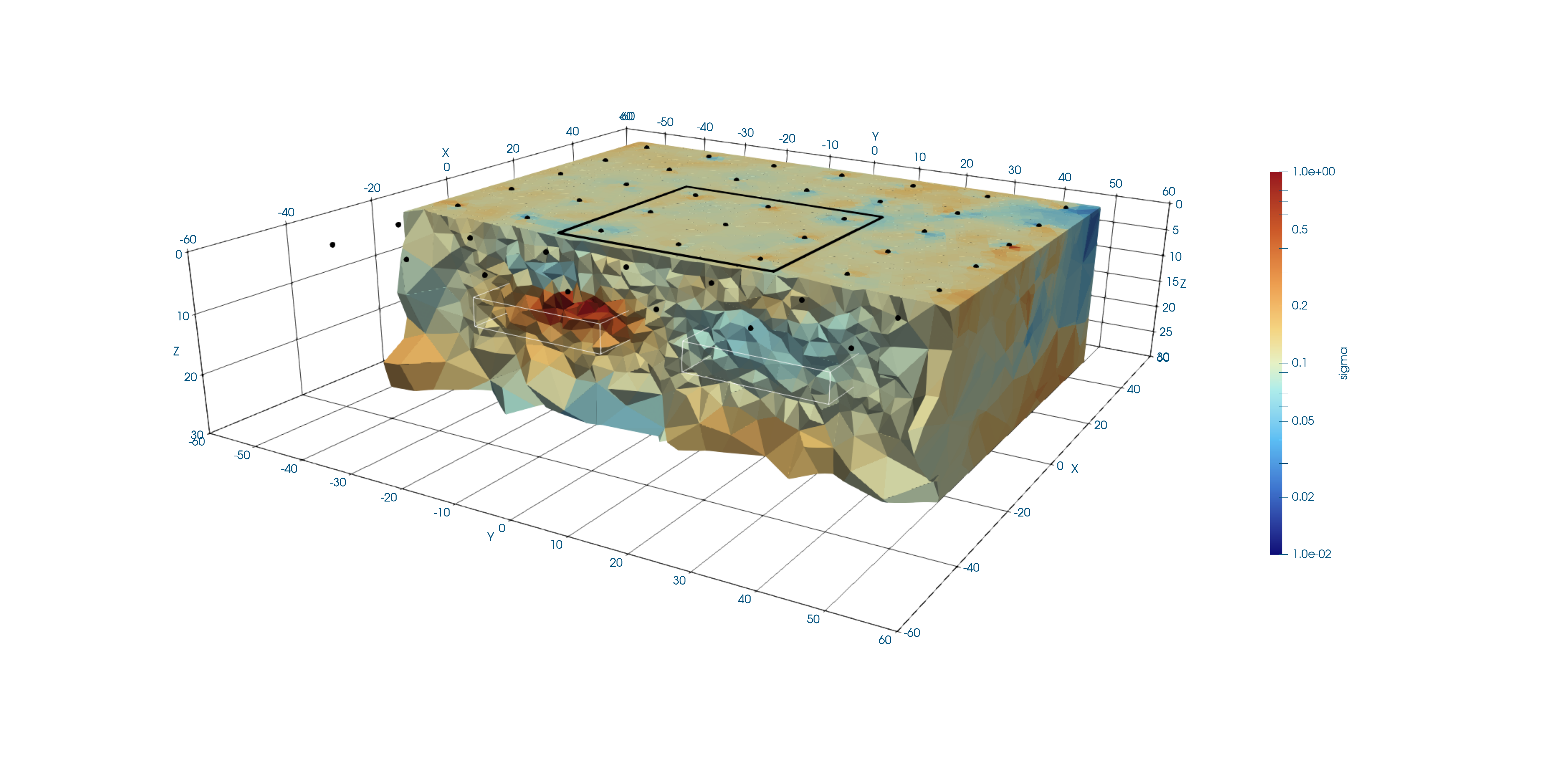}
    \vspace*{-12mm}
    \caption{Slice through the reconstructed conductivity model obtained after 25 Gauss--Newton iterations. The slice illustrates the lateral extent and depth positioning of the recovered anomalies. While the main structures are clearly identified, the anomalies are moderately smoothed and exhibit reduced contrast, reflecting the influence of regularization and the limited resolution of the data.}
    \label{fig:model-result-sliced}
\end{figure}

The reconstructed model obtained from the synthetic data is shown in Fig. \ref{fig:model-result} (horizontal slice) and Fig. \ref{fig:model-result-sliced} (vertical slice).
The inversion recovers the principal features of the true conductivity distribution.
In particular, the locations of the conductive blocks are correctly identified, and their lateral extent is reproduced with reasonable accuracy.

The recovered anomalies exhibit some degree of smoothing compared to the true model. 
This behaviour is expected due to the applied smoothness regularization and the limited resolving power of the TEM data, especially at depth. 

The conductive anomalies are well recovered in both geometry and magnitude, whereas the resistive anomalies exhibit a reduced contrast and are only weakly resolved. 
This behaviour is consistent with the limited sensitivity of electromagnetic data to resistive structures and the regularization-induced smoothing, which together lead to an underestimation of resistive amplitudes.

In the vertical direction, the anomalies indicate limited depth resolution. 
This effect is typical for near-surface TEM configurations, where early-time data constrain shallow structures more strongly than deeper parts of the model. 
Nevertheless, the inversion places the anomalies at approximately the correct depth.

The spatial distribution of artefacts is limited, and no significant spurious structures are observed away from the target regions. 
This indicates that the combination of data weighting and regularization provides a stable inversion result.

Overall, the example demonstrates that the proposed inversion framework is capable of recovering the dominant features of a 3-D conductivity distribution from transient electromagnetic data, while exhibiting the expected characteristics associated with regularized inverse problems.

\begin{figure}
    \centering    \includegraphics[width=0.5\linewidth]{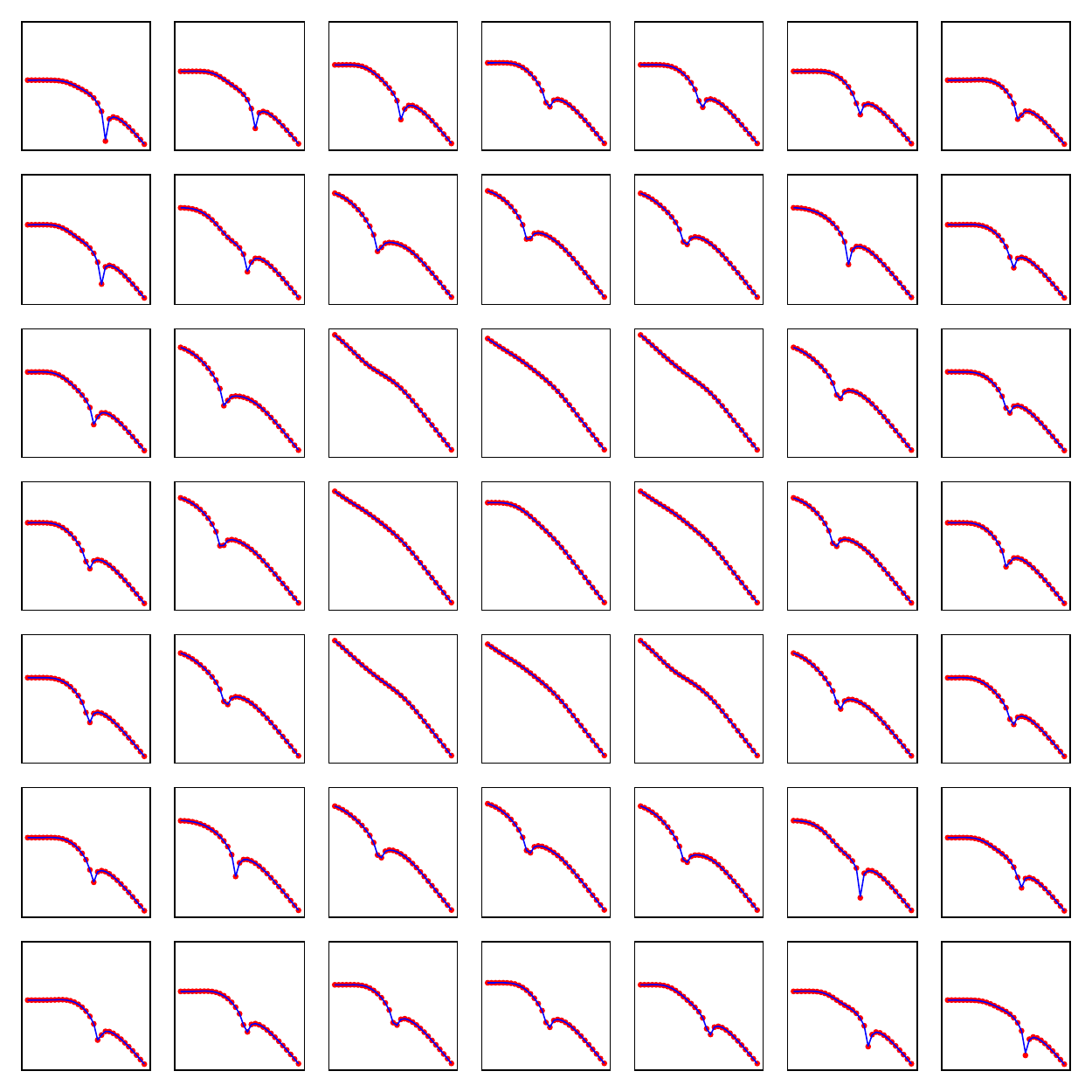}
    
    \caption{Comparison of observed and predicted TEM responses ($\partial_t B_z)$ after 25 iterations for all receiver locations, arranged according to the $7 \times 7$ acquisition grid. 
    Each panel shows one transient, with observed data indicated by red symbols and the model response by blue lines. 
    All panels share identical axes, with time ranging from $10^{-6}$ to $10^{-3}$ s and amplitudes from $10^{-8}$ to $10^{-2}$ $\mathrm{V\,A^{-1}\,m^{-2}}$. 
    Axes and annotations are omitted for clarity, emphasizing the overall agreement in shape and amplitude across the survey area.}
    \label{fig:fit-plot}
\end{figure}

To provide a compact overview of the data fit across all receiver locations, we employ a simplified sparkline visualization in which the observed and predicted transients are displayed as small multiples arranged according to the spatial layout of the receiver grid. All graphical decorations are omitted to emphasize the overall agreement in response shape and amplitude. This representation allows for a rapid qualitative assessment of the inversion result over the entire survey area (Fig.~\ref{fig:fit-plot}).

\begin{figure}
    \centering
    \includegraphics[width=0.7\linewidth]{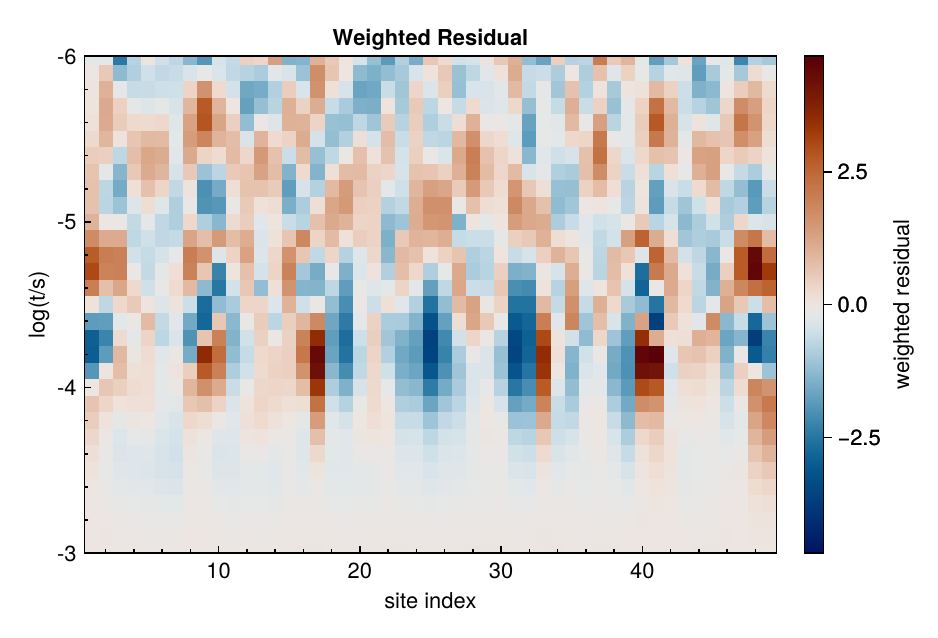}
    \caption{Heat map of the weighted data residuals $\mathbf{W}_d(\mathbf{d}(\mathbf{m}) - \mathbf{d}^{\mathrm{obs}})$ after 25 Gauss--Newton iterations. 
    Rows correspond to time channels, and columns correspond to receiver locations arranged according to the acquisition grid. 
    The colour scale is symmetric about zero, highlighting positive and negative weighted deviations of the model response from the observed data.}
    \label{fig:weighted-residual}
\end{figure}

The weighted data residuals $\mathbf{W}_d(\mathbf{d}(\mathbf{m}) - \mathbf{d}^{\mathrm{obs}})$ after 25 Gauss--Newton iterations across all receivers and time channels are displayed in Fig.~\ref{fig:weighted-residual}.
The distribution of the misfit in this representation exhibits a random pattern which is in good agreement with the considered error model.  

\begin{figure}
    \centering
    \includegraphics[width=0.99\linewidth]{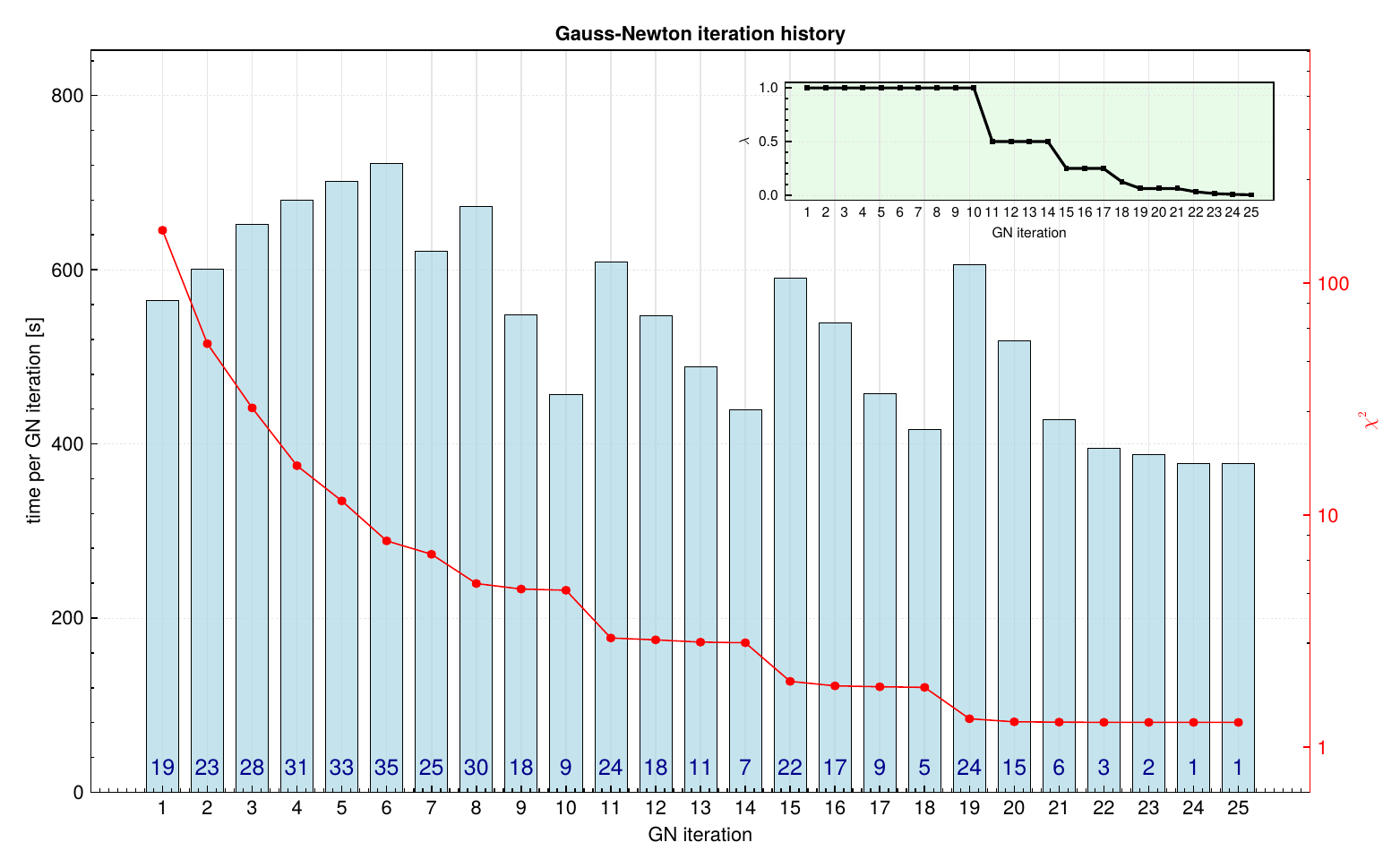}
    \caption{Inversion convergence and performance. Bars show the computational time required for each Gauss–Newton (GN) iteration (left axis), while the red curve indicates the corresponding evolution of the data misfit $\chi^2$ (right axis). Numerical annotations within the bars denote the number of LSQR iterations performed at each Gauss–Newton step. The inset displays the evolution of the regularization parameter $\lambda$ throughout the inversion. 
    }
    \label{fig:fig-chi2}
\end{figure}

The computational performance is summarized in Fig.~\ref{fig:fig-chi2}, which shows both the runtime per Gauss–Newton iteration and the corresponding evolution of the data misfit $\chi^2$. 
The figure demonstrates that the computational cost is dominated by the factorization and repeated Jacobian-related operations required within the inner LSQR iterations.

A Gauss–Newton iteration with one inner LSQR iteration required approximately 380 s, while the most expensive step took about 720 s and involved 35 LSQR iterations. 
These observations will be discussed in more detail in the next section.
The peak memory consumption was approximately 620 GB, consistent with the storage of all 21 factorized shifted systems in memory.

\subsection{Computational performance}

To assess the parallel performance of the proposed RBA formulation, we conducted a series of scaling experiments using 16 poles in the rational approximation. 
This choice is computationally convenient because 16 admits balanced decompositions across 2, 4, 8, and 16 MPI ranks, thereby allowing the pole-wise parallelism of the method to be examined systematically.
In the inversion experiments, however, 21 poles were employed in order to obtain a more accurate approximation of the transient response over the considered time interval. 
The scaling experiments should therefore primarily be interpreted as a controlled illustration of the parallel characteristics of the method, whereas the inversion timings reported below correspond to the more accurate 21-pole configuration.

\begin{figure}
    \centering
    \includegraphics[width=0.9\linewidth]{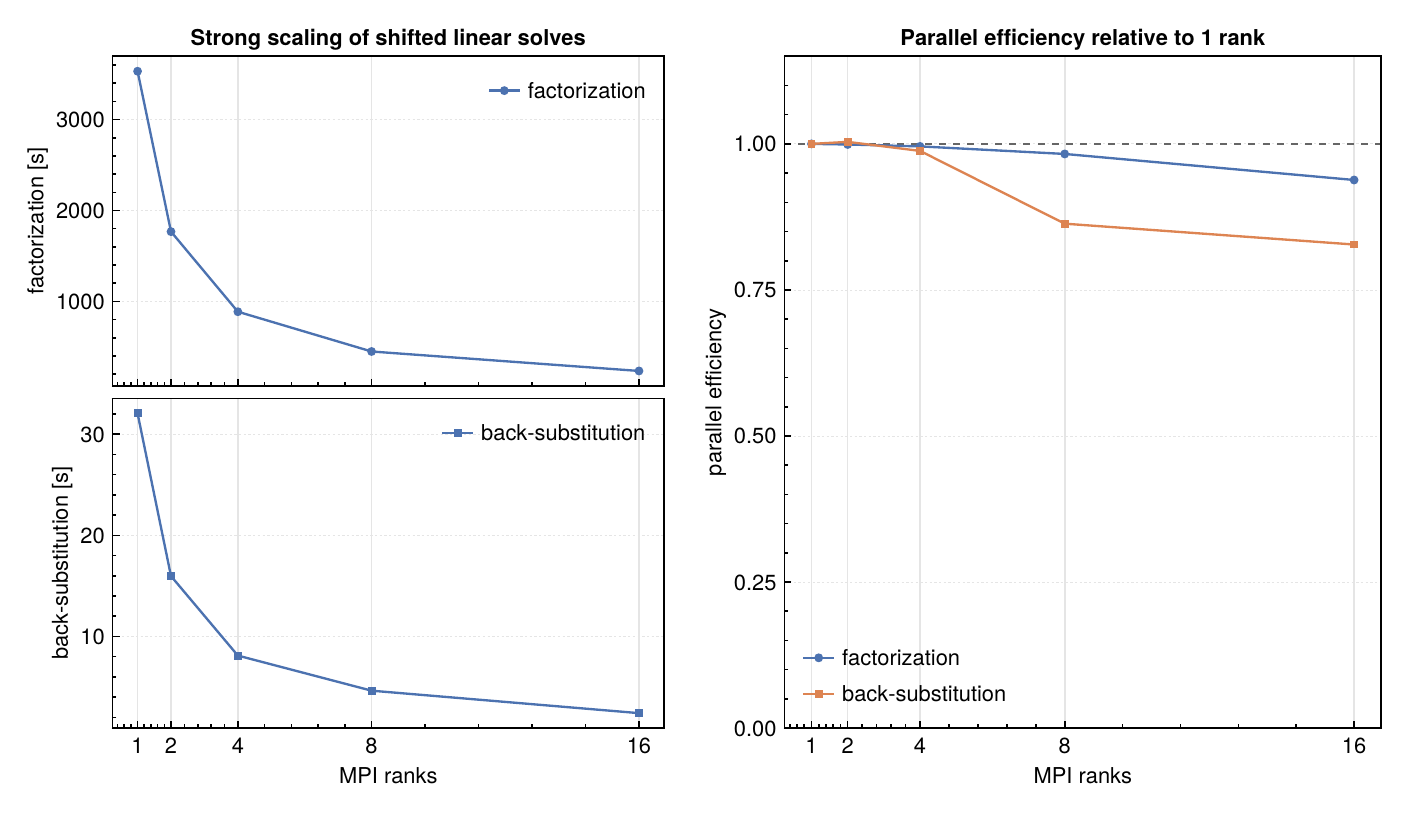}
    \caption{Strong-scaling behaviour of the shifted linear systems for an RBA approximation with 16 poles. The left panels show the measured runtimes for the sparse factorization and subsequent back-substitution stages, while the right panel shows the corresponding parallel efficiency relative to the single-rank run.}
    \label{fig:par-efficiency}
\end{figure}

The results illustrated in Fig.~\ref{fig:par-efficiency} demonstrate near-ideal strong scaling up to eight MPI ranks for both the factorization and the subsequent triangular solves, confirming that the shifted systems can be processed largely independently. 
A moderate degradation in parallel efficiency is observed for the largest run with 16 ranks. 
Since all shifted systems possess identical sparsity structure, this loss of efficiency is not related to changes in fill-in or elimination-tree complexity, but is primarily attributable to the shared-memory NUMA architecture used in the experiments. 
In the current implementation, the MUMPS factors are stored in shared memory and accessed concurrently by multiple MPI ranks. 
As the number of ranks increases, contention for memory bandwidth and cache resources becomes increasingly significant, particularly during the sparse triangular solves, which are comparatively memory-bandwidth dominated and exhibit lower arithmetic intensity than the factorization phase. 
Consequently, the reduction in local computational work per rank is eventually offset by increased memory-access overhead. 
This limitation is expected to be less pronounced on distributed-memory HPC systems, where the shifted systems can be distributed across independent nodes and therefore interfere less strongly within the memory subsystem.

\begin{figure}
    \centering
    \includegraphics[width=0.7\linewidth]{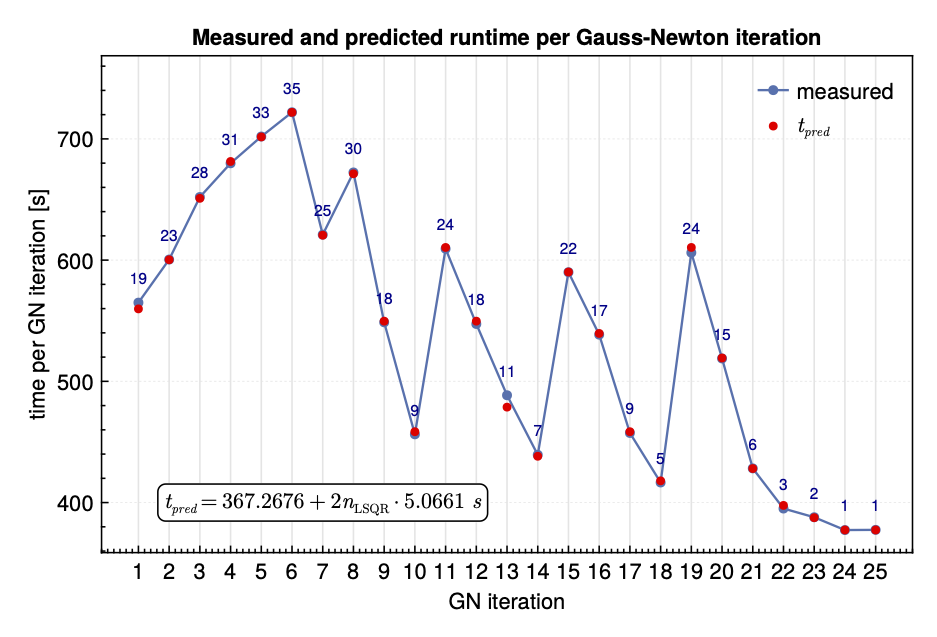}
    \caption{Measured and predicted runtimes per Gauss–Newton iteration for an RBA approximation with 21 poles. The predicted runtimes are obtained from a linear performance model consisting of a fixed factorization cost and a contribution proportional to the number of LSQR iterations. Numerical labels indicate the corresponding number of LSQR iterations for each Gauss–Newton step.}
    \label{fig:pred-runtime}
\end{figure}

The same performance characteristics extend naturally to the inversion stage. 
The measured runtimes of the Gauss–Newton iterations are accurately described by a simple linear performance model consisting of a fixed factorization cost plus a contribution proportional to the number of LSQR iterations (Fig.~\ref{fig:pred-runtime}). 

This confirms that the computational effort of the inversion is dominated by the repeated triangular solves associated with the evaluation of $\mathbf J\vv$ and $\mathbf J^{\top}\mathbf w$. 
Since these operations inherit the same pole-wise parallelism as the forward problem, the observed scalability of the shifted solves directly translates into reduced wall-clock times for the complete inversion workflow.

\subsection{Comparison to other numerical approaches}

To place the computational performance of the proposed RBA formulation into context, it is instructive to compare its numerical complexity with that of conventional implicit time-stepping approaches commonly used in time-domain EM inversion. 
In a previous study \citep{rochlitz2021evaluation}, three different transient simulation strategies were investigated, namely implicit Euler time stepping, Fourier-based synthesis, and a Krylov-subspace method based on a rational Arnoldi approximation of the matrix exponential. 
Among these approaches, implicit Euler remains the most relevant reference for inversion applications because of its widespread use in existing inversion frameworks. 
For a representative transient with 31 time gates between $10^{-6}$ and $10^{-3}\,\mathrm{s}$, the implicit Euler formulation requires 31 real-valued sparse factorizations, each followed by approximately 10 inexpensive time steps \citep{rochlitz2021evaluation}. 
Thus, although the factorizations can be reused locally for several time steps, the computation remains fundamentally tied to a sequential time-marching procedure.

In the inversion experiments presented here, 42 poles were used in the rational approximation in order to obtain higher accuracy over the considered transient interval. 
Exploiting conjugate symmetry reduces the corresponding computational problem from 42 to 21 complex-valued shifted systems. 
Compared with the implicit Euler approach of \cite{rochlitz2021evaluation}, which requires 31 real-valued factorizations for the considered transient discretization, the total serial factorization effort of the RBA formulation is therefore of comparable magnitude even when the additional cost of complex arithmetic is taken into account. 
However, the decisive difference between the two approaches lies not in the nominal number of factorizations, but in the computational structure of the resulting algorithms. 
Implicit Euler remains fundamentally tied to a sequential time-stepping procedure, whereas all shifted systems of the RBA formulation are mutually independent and can therefore be processed concurrently. 
Moreover, the computational cost of the RBA approach is independent of the number of requested time channels and avoids the sequential forward/backward recursions required for the evaluation of $\mathbf J \vv$ and $\mJ^{\top}\mathbf w$ in implicit time-stepping schemes. 
Consequently, the principal advantage of the RBA formulation is the substantial reduction of the parallel critical path, which leads to significantly improved wall-clock performance on modern distributed-memory architectures.

\section{Conclusion}

In this study, we presented a parallel implementation for large-scale 3-D electromagnetic inversion and demonstrated its capability to successfully recover the conductivity structure of a synthetic 3-D model problem. 
The proposed approach combines a classical Gauss–Newton inversion framework with the Rational Best-Approximation and parallel sparse direct solvers. 
Parallel scalability was demonstrated on a shared-memory architecture, showing that the computational workload associated with the shifted systems can be efficiently distributed across multiple MPI processes.

The implementation was realized using a finite-element framework; however, the proposed inversion strategy itself is independent of any specific finite-element library and can therefore be transferred to other discretization environments. Model regularization was formulated using Raviart-Thomas $RT_0$ elements, which provided smooth and physically consistent conductivity reconstructions. 
Furthermore, LSQR was employed for the solution of the linearized inverse problems and proved advantageous compared to conventional conjugate-gradient approaches in terms of convergence behaviour and robustness.

The numerical experiments also revealed the principal computational bottlenecks of the current implementation. 
In particular, the simultaneous storage of multiple sparse matrix factorizations leads to a substantial memory footprint, while the execution of many parallel processes on a shared-memory machine introduces additional limitations related to memory bandwidth and resource contention. 
Future improvements will  focus on replacing direct factorizations with suitably preconditioned iterative solvers in order to reduce memory consumption and improve scalability.

Although the RBA formulation uses complex-valued shifted systems, which are more expensive than the real systems arising in standard implicit time stepping, this does not remove the main advantage of the method. 
A complex sparse solve may be estimated as costing roughly two to four real solves, depending on implementation and hardware. 
Even under this conservative assumption, the RBA approach replaces a sequential recursion over all time steps by a fixed number of mutually independent shifted solves. 
Therefore, the principal gain is in the reduced parallel critical path and in the independence of the cost from the number of time channels, rather than in a simple reduction of serial arithmetic work.

An important next step will be the extension of the implementation towards distributed high-performance computing clusters, enabling significantly larger inverse problems and improved parallel efficiency. 
Finally, the release of the software as open-source code (URL will be made public upon publication) is a key aspect of this work, promoting reproducibility, transparency and further methodological development within the community.

\bibliography{references}

\end{document}